\documentclass[12pt]{amsart}
\usepackage{epsfig,color}

\makeatother

\theoremstyle{definition}

\def\fnum{equation}
\newtheorem{Thm}[\fnum]{Theorem}

\newtheorem{Lem}[\fnum]{Lemma}

\numberwithin{equation}{section}

\newcommand{\dist}{{\text {dist}}}

\def\RR{{\bold R}}
\def\SS{{\bold S}}

\newcommand{\Length}{{\text {Length}}}
\newcommand{\Energy}{{\text {Energy}}}

\newcommand{\eqr}[1]{(\ref{#1})}

\begin{document}

\title[Min-max for sweepouts by curves]{Min-max for sweepouts by curves}

\author{Tobias H. Colding}%
\address{MIT\\
77 Massachusetts Avenue, Cambridge, MA 02139-4307\\
and Courant Institute of Mathematical Sciences\\
251 Mercer Street, New York, NY 10012.}
\author{William P. Minicozzi II}%
\address{Department of Mathematics\\
Johns Hopkins University\\
3400 N. Charles St.\\
Baltimore, MD 21218}

\thanks{The   authors
were partially supported by NSF Grants DMS  0606629 and DMS
0405695}


\email{colding@math.mit.edu  and minicozz@math.jhu.edu}

\maketitle

\section{Introduction}

 Given a Riemannian metric on the $2$-sphere, sweep the
$2$-sphere out  by a continuous one-parameter family of closed
curves starting and ending at point curves. Pull the sweepout
tight by, in a continuous way, pulling each curve as tight as
possible yet preserving the sweepout.  We show the following
useful property (see Theorem \ref{t:mm} below); cf. \cite{CM1},
\cite{CM2}, proposition 3.1 of \cite{CD}, proposition 3.1 of
\cite{Pi}, and 12.5 of \cite{Al}:

\vskip2mm
\parbox{6in}{Each curve in the tightened
sweepout whose length is close to the length of the longest curve
in the sweepout must itself be close to a closed geodesic. In
particular, there are curves in the sweepout that are close to
closed geodesics.}

\vskip2mm \noindent Finding closed geodesics on the $2$-sphere by
using sweepouts goes back to Birkhoff in the 1920s; see \cite{B}
and  section $2$ in \cite{Cr} about  Birkhoff's ideas. The
argument works equally well on any closed manifold, but only
produces non-trivial closed geodesics when the width, which is
defined in \eqr{e:w} below, is positive.  For instance, when $M$
is topologically a $2$-sphere, the width is loosely speaking the
length of the shortest closed curve needed to ``pull over'' $M$.
Thus Birkhoff's argument gives that the width is realized as the
length of a closed geodesic.

The above useful property is virtually always implicit in any
sweepout construction of critical points for variational problems
yet it is not always recorded since most authors are only
interested in the existence of one critical point.

Similar results holds for sweepouts by $2$-spheres instead of
circles; cf. \cite{CM2}.  The ideas are essentially the same in
the two cases, though the techniques in the curve case are purely
ad hoc whereas in the $2$-sphere case  additional techniques,
developed in the 1980s, have to be used to deal with   energy
concentration (i.e., ``bubbling''); cf. \cite{Jo}.

\section{Existence of good sweepouts by curves}
\label{s:existe}

Let $M$ be a closed Riemannian manifold.
 Fix a large positive integer $L$ and let $\Lambda$ denote
the space of piecewise linear maps from $\SS^1$ to $M$ with
exactly $L$ breaks (possibly with unnecessary breaks) such that
the length of each geodesic segment is at most $2\pi$,
parametrized by a (constant) multiple of arclength, and with
Lipschitz bound $L$. By a linear map, we mean a (constant speed)
geodesic. Let $G \subset \Lambda$ denote the set of immersed
closed geodesics in $M$ of length at most $2\pi L$. (The energy of
a curve in $\Lambda$ is equal to its length squared divided by
$2\pi$. In other words, energy and length are essentially
equivalent.)

We will use the distance and topology on $\Lambda$ given by the
$W^{1,2}$ norm (Sobolev norm) on the space of maps from $\SS^1$ to
$M$. The simplest way to define the $W^{1,2}$ norm is to
isometrically embed the compact manifold $M$ into some Euclidean
space $\RR^N$.{\footnote{Recall that the square of the $W^{1,2}$
norm of a map $f:\SS^1 \to \RR^N$ is
\begin{equation}
     \int_{\SS^1}  \left( |f|^2 + |f'|^2
    \right)  \, . \notag
\end{equation}
Thus two curves that are $W^{1,2}$ close are also $C^0$ close; cf.
\eqr{e:holder}.}} It will be convenient to scale $\RR^N$, and thus
$M$, by a constant so that it satisfies the following: (M1)
$\sup_M |A| \leq 1/16$, where $|A|^2$ is the norm squared of the
second fundamental form of $M$, i.e., the sum of the squares of
the principal curvatures (see, e.g., (1.24) on page 4 of
\cite{CM3}); (M2)  the injectivity radius of $M$ is at least
$8\pi$ and the curvature is at most $1/64$, so that every geodesic
ball of radius at most $4\pi$ in $M$ is strictly geodesically
convex; (M3) if $x,y \in M$ with $|x-y|\leq 1$, then
  $\dist_M (x,y) \leq 2|x-y|$.

\subsection{The width}     \label{ss:width}

Let $\Omega$ be the set of continuous maps $\sigma : \SS^1 \times [-1,1] \to M$ so that
 for each $t$ the map $\sigma (\cdot , t )$
is in $W^{1,2}$, the map  $t \to \sigma (\cdot , t )$
is continuous
 from $[-1,1]$ to
$W^{1,2}$, and finally
$\sigma$ maps $\SS^1 \times \{ -1 \}$ and $\SS^1 \times \{ 1 \}$ to points.
Given a map
    $\hat{\sigma} \in \Omega$, the homotopy class $\Omega_{\hat{\sigma}}$
 is defined to be the set of maps $\sigma \in \Omega$ that are homotopic to
    $\hat{\sigma}$ through maps in $\Omega$.
The width $W = W (\hat{\sigma})$ associated to the homotopy class
 $\Omega_{\hat{\sigma}}$ is defined by taking $\inf$ of $\max$ of
the energy of each slice.{\footnote{A particularly interesting
example is when $M$ is a topological $2$-sphere and the induced
map from $\SS^2$ to $M$ has degree one.  In this case, the width,
defined below, is positive and realized by one or more non-trivial
closed geodesics. In general, the width is always non-negative but
may not always be  positive.}}
That is,  set
\begin{equation}    \label{e:w}
    W = \inf_{  \sigma \in \Omega_{\hat{\sigma}}  } \,
       \,  \max_{ t \in [-1, 1]} \,  \Energy \, (\sigma (\cdot , t ))
          \, ,
\end{equation}
where the energy is given by $\Energy \, (\sigma (\cdot , t )) =
\int_{\SS^1} \, \left| \partial_x \sigma (x,t) \right|^2 \, dx$.

The main theorem, Theorem \ref{t:mm},  that almost maximal slices in the
tightened sweepout are almost geodesics, is proven in subsection \ref{ss:mm}.  The proof of
this theorem as well as the
construction of the sequence of tighter and tighter sweepouts
uses a curve shortening map
that is defined in the next subsection.  We also state the key properties of
the shortening map in the next subsection, but postpone their proofs
to Section \ref{s:A} and the appendices.

\subsection{Curve shortening $\Psi$}    \label{ss:birk}

The curve shortening is a map $\Psi: \Lambda \to \Lambda$ so
that{\footnote{This map is essentially what is usually called
Birkhoff's curve shortening process, see section 2 of \cite{Cr}.}}
\begin{itemize}
\item[(1)] $\Psi(\gamma)$ is
homotopic to $\gamma$ and $\Length (\Psi(\gamma)) \leq \Length
(\gamma)$.
\item[(2)] $\Psi (\gamma)$ depends continuously on $\gamma$.
\item[(3)] There is a continuous function $\phi:[0,\infty) \to
[0,\infty)$ with $\phi (0) = 0$ so that
\begin{equation}    \label{e:01}
    \dist^2 (\gamma , \Psi (\gamma)) \leq \phi \left( \frac{\Length^2 (\gamma) -
\Length^2 (\Psi (\gamma))}{\Length^2 (\Psi (\gamma))} \right) \, .
\end{equation}
\item[(4)] Given $\epsilon > 0$, there exists $\delta > 0$ so that
if $\gamma \in \Lambda$ with $\dist (\gamma , G) \geq \epsilon$,
then $\Length \, (\Psi (\gamma)) \leq \Length \, (\gamma) -
\delta$.
\end{itemize}

  To define $\Psi$, we will fix a partition of   $\SS^1$
by choosing $2L$ consecutive evenly spaced points{\footnote{Note
that this is not necessarily where the piecewise linear maps have
breaks.}}
\begin{equation}
    x_0 , x_1 , x_2, \dots , x_{2L} = x_0 \in \SS^1 \, ,
\end{equation}
so that $|x_{j} - x_{j+1}| = \frac{\pi}{L}$.  $\Psi (\gamma)$ is
given in three steps.  First, we apply step 1 to $\gamma$ to get a
curve $\gamma_e$, then we apply step 2 to $\gamma_e$ to get a
curve $\gamma_o$.  In the third and final step, we reparametrize
$\gamma_o$ to get $\Psi (\gamma)$.

\vskip1mm \noindent {\bf{Step 1}}: Replace $\gamma$ on each
{\emph{even}} interval, i.e., $[x_{2j} , x_{2j+2}]$, by the linear
map with the same endpoints to get a piecewise linear curve
$\gamma_e: \SS^1 \to M$. Namely, for each $j$, we let $\gamma_e
\big|_{[x_{2j},x_{2j+2}]}$ be the unique shortest (constant speed)
geodesic from $\gamma(x_{2j})$ to $\gamma(x_{2j+2})$.

\vskip1mm \noindent {\bf{Step 2}}: Replace $\gamma_e$ on each
{\emph{odd}} interval by the linear map with the same endpoints to
get the piecewise linear curve $\gamma_o: \SS^1 \to M$.

\vskip1mm \noindent {\bf{Step 3}}: Reparametrize $\gamma_o$
(fixing $\gamma_o (x_0)$) to get the desired  constant speed
curve $\Psi(\gamma) : \SS^1 \to M$.

It is easy to see that $\Psi$ maps $\Lambda$ to $\Lambda$ and has
property (1); cf. section 2 of \cite{Cr}. Properties (2), (3) and
(4) for $\Psi$ are established in Section \ref{s:A} and Appendix
\ref{s:aB}.  Throughout the rest of this section, we will assume
these properties and use them to prove the main theorem.

\vskip2mm  The next lemma, which combines  (3) and (4), is the key
to producing the desired sequence of sweepouts.

\begin{Lem} \label{l:i1}
Given $W \geq 0$ and $\epsilon > 0$, there exists $\delta > 0$  so that
if $\gamma \in \Lambda$ and
\begin{equation}    \label{e:closei}
   2\pi \, (W - \delta) <  \Length^2 \, ( \Psi (\gamma) ) \leq  \Length^2 \, ( \gamma)  < 2\pi \, (W + \delta ) \, ,
\end{equation}
then  $\dist ( \Psi (\gamma) , G) < \epsilon$.
\end{Lem}

\begin{proof}
If $W \leq \epsilon^2/6$, then the Wirtinger inequality (see footnote $6$) yields the lemma with $\delta = \epsilon^2/6$.

Assume next that $W > \epsilon^2/6$.
 The triangle
inequality    gives
\begin{equation}    \label{e:closei2}
    \dist ( \Psi (\gamma) , G) \leq
        \dist ( \Psi (\gamma) , \gamma) + \dist (  \gamma , G)
     \, .
\end{equation}
Since $\Psi$  does not decrease  the length of $\gamma$ by much,
property (4) of $\Psi$ allows us to bound $\dist ( \gamma , G)$ by
$\epsilon/2$ as long as $\delta$ is sufficiently small. Similarly,
property (3) of $\Psi$ allows us to bound $\dist ( \Psi (\gamma) ,
\gamma)$ by $\epsilon/2$ as long as $\delta$ is sufficiently
small.
\end{proof}

\subsection{Defining the sweepouts}     \label{ss:sweep}

Choose a sequence of maps $\hat{\sigma}^j \in
\Omega_{\hat{\sigma}}$
 with
\begin{equation}    \label{e:minseq}
    \max_{t \in [-1,1]} \, \, \Energy \, (\hat{\sigma}^j (\cdot , t))
    < W + \frac{1}{j} \, .
\end{equation}
Observe that \eqr{e:minseq} and the Cauchy-Schwarz inequality
imply a uniform bound for the length and uniform $C^{1/2}$
continuity for the slices,   that are both independent of   $t$
and $j$. The first follows immediately and the latter follows from
\begin{align}   \label{e:holder}
      \left| \hat{\sigma}^j (x , t) \, - \,  \hat{\sigma}^j (y ,
    t)\right|^2 &\leq \left( \int_x^y \left| \partial_s \hat{\sigma}^j (s , t)
    \right| \, ds \right)^2 \notag \\
    &\leq |y-x| \,   \int_x^y \left| \partial_s \hat{\sigma}^j (s , t)
    \right|^2 \, ds  \leq |y-x| \, (W+1) \, .
\end{align}

 We will replace the $\hat{\sigma}^j$'s by sweepouts
$\sigma^j$ that, in addition to  satisfying \eqr{e:minseq}, also
satisfy that the slices $\sigma^j (\cdot , t)$ are in $\Lambda$.
  We will do this by using   local
linear replacement similar to   Step 1 of the construction of
$\Psi$.  Namely, the uniform $C^{1/2}$ bound for the slices allows
us to fix a partition of points $y_0 , \dots , y_N = y_0$ in
$\SS^1$ so that each interval $[y_i , y_{i+1}]$ is always mapped
to a ball in $M$ of radius at most  $4\pi$.  Next, for each $t$ and each $j$, we
replace $\hat{\sigma}^j (\cdot , t) \, \big|_{[y_{i},y_{i+1}]}$ by
the linear map (geodesic) with the same endpoints  and call the
resulting map $\tilde{\sigma}^j (\cdot , t)$.  Reparametrize $\tilde{\sigma}^j (\cdot , t)$ to have constant
speed to get $\sigma^j (\cdot , t)$.
 It is easy to see that each $\sigma^j (\cdot , t)$
satisfies \eqr{e:minseq}. Furthermore, the length bound for
$\sigma^j (\cdot , t)$ also gives a uniform Lipshitz bound for the
linear maps; let $L$ be the maximum of $N$ and this Lipshitz bound.

  It remains to show that $\sigma^j$ is continuous in the transversal
direction, i.e., with respect to $t$, and homotopic to $
\hat{\sigma}$ in $\Omega$.
 These facts were established
already by Birkhoff (see \cite{B} and section $2$ of \cite{Cr}),
but also follow immediately from Appendix \ref{s:aB}.

Finally, applying the  replacement map $\Psi$ to   each $\sigma^j
(\cdot , t)$ gives a new  sequence of sweepouts $ \gamma^j = \Psi
(\sigma^j)$.  (By Appendix \ref{s:aB},  $\Psi$ depends
continuously on $t$ and preserves the homotopy class
$\Omega_{\hat{\sigma}}$; it is clear that $\Psi$
 fixes the constant maps at $t = \pm 1$.)

\subsection{Almost maximal implies almost critical}
\label{ss:mm}

Our main result is that this sequence $\gamma^j$ of sweepouts is tight in the
sense of the Introduction.  Namely, we have the following
theorem.

\begin{Thm}     \label{t:mm}
Given $W \geq 0$ and $\epsilon > 0$, there exist $\delta > 0$   so that
if $j
> 1/\delta$ and for some $t_0$
\begin{equation}    \label{e:close1}
    2\pi \, \Energy \, ( \gamma^j (\cdot , t_0))
        = \Length^2 \, ( \gamma^j (\cdot , t_0)) > 2\pi \, (W - \delta)  \, ,
\end{equation}
then for this $j$ we have $\dist \, \left(  \gamma^j (\cdot , t_0)
\, , \, G \right) < \epsilon$.
\end{Thm}

 \begin{proof}
 Let $\delta$ be given by
Lemma \ref{l:i1}.  By \eqr{e:close1}, \eqr{e:minseq}, and using
that $j
> 1/\delta$, we get
\begin{equation}    \label{e:close1a}
2\pi \, (W - \delta) <  \Length^2 \, ( \gamma^j (\cdot , t_0))
\leq    \Length^2 \, ( \sigma^j (\cdot , t_0)) <   2\pi \, (W+
\delta)  \, .
\end{equation}
Thus, since $\gamma^j (\cdot , t_0) = \Psi ( \sigma^j (\cdot ,
t_0))$, Lemma \ref{l:i1}   gives $\dist ( \gamma^j (\cdot , t_0)
\, , \, G) < \epsilon$, as claimed.
\end{proof}

\section{Establishing Properties (2), (3) and (4) for $\Psi$}  \label{s:A}

To prove  (2) and (3), it is useful to observe that there is an
equivalent, but more symmetric, way to construct $\Psi (\gamma)$
using four steps:
\begin{enumerate}
\item[($A_1$)] Follow Step 1 to get $\gamma_e$. \item[($B_1$)]
Reparametrize $\gamma_e$ (fixing the image of $x_0$) to get the
constant speed curve $\tilde{\gamma}_e$.  This reparametrization
moves the points $x_j$ to new points $\tilde{x}_j$ (i.e.,
$\gamma_e (x_j) = \tilde{\gamma}_e (\tilde{x}_j)$). \item[($A_2$)]
Do linear replacement on the odd $\tilde{x}_j$ intervals to get
$\tilde{\gamma}_o$. \item[($B_2$)] Reparametrize
$\tilde{\gamma}_o$  (fixing the image of $x_0$) to get the
constant speed curve $\Psi (\gamma)$.
\end{enumerate}
The reason that this gives the same curve is that
$\tilde{\gamma}_o$ is just a reparametrization of  ${\gamma}_o$.
We will also use that each of the four steps is energy
non-increasing.  This is obvious for  the linear replacements,
since linear maps minimize energy.    It follows  from the
Cauchy-Schwarz inequality for the reparametrizations, since for a
curve $\sigma:\SS^1 \to M$ we have
\begin{equation}
    \Length^2 (\sigma) \leq 2\pi \, \Energy (\sigma) \, ,
\end{equation}
with equality if and only if $|\sigma'| = \Length (\sigma)/(2\pi)$
almost everywhere.

\vskip1mm Using the alternative way of defining $\Psi (\gamma)$ in
four steps, we see that (3) follows from the triangle inequality
once   we bound $\dist (\gamma , \gamma_e)$ and $\dist (\gamma_e ,
\tilde{\gamma}_e)$ in terms of the decrease in length  (as well as
the analogs for steps $(A_2)$ and $(B_2)$).

The bound on $\dist (\gamma , \gamma_e)$  follows directly from
the following, see Appendix \ref{s:B} for the proof:

\begin{Lem} \label{l:sc}
There exists $C$ so that if  $I$ is an interval of length at most
$2\pi/L$, $\sigma_1: I\to M$ is a Lipschitz curve with
$|\sigma_1'| \leq L$, and $\sigma_2 : I\to M$ is the minimizing
geodesic with the same endpoints, then
\begin{equation}    \label{e:a2}
    \dist^2 (\sigma_1 , \sigma_2) \leq C \, \left( \Energy (\sigma_1) -
    \Energy (\sigma_2) \right) \, .
\end{equation}
\end{Lem}

Applying Lemma \ref{l:sc} on each of the $L$ intervals  in step
$(A_1)$, we get that
\begin{equation}    \label{e:a2vv}
    \dist^2 (\gamma , \gamma_e) \leq C \, \left( \Energy (\gamma) -
    \Energy (\gamma_e) \right) \leq \frac{C}{2\pi} \, \left( \Length^2 (\gamma) -
    \Length^2 (\Psi(\gamma)) \right) \, \, .
\end{equation}
 This gives the desired bound on $\dist (\gamma , \gamma_e)$ since
 $\Length (\Psi (\gamma)) \leq 2\pi \, L$.

\vskip1mm
 In bounding $\dist (\gamma_e , \tilde{\gamma}_e)$, we will use that
  $\gamma_e$ is just the
composition $\tilde{\gamma}_e \circ P$, where $P: \SS^1 \to \SS^1$
is a monotone piecewise linear map.{\footnote{The map $P$ is
Lipschitz, but the inverse map $P^{-1}$ may not be if $\gamma_e$
is constant on an interval.}} Using that $|\tilde{\gamma}_e'| =
\Length (\tilde{\gamma}_e)/(2\pi)$ (away from the breaks) and that
the integral of $P'$ is $2\pi$, an easy calculation gives
\begin{align}   \label{e:ta}
    \int \left( P' - 1 \right)^2 &= \int ( P')^2 - 2\pi =
\int \left( \frac{|\gamma_e'|}{|\tilde{\gamma}_e' \circ P|}
\right)^2 - 2\pi    =
      \frac{4\pi^2 }{\Length^2
(\tilde{\gamma}_e)}    \, \int |\gamma_e'|^2 - 2\pi \notag
\\ &= 2\pi \,   \frac{ \Energy
    (\gamma_e) - \Energy (\tilde{\gamma}_e)}{\Energy
    (\tilde{\gamma}_e)}    \leq
2\pi \,   \frac{ \Energy
    (\gamma) - \Energy (\Psi (\gamma))}{\Energy
    (\Psi (\gamma))}
    \, .
\end{align}
Since $\gamma_e$ and $\tilde{\gamma}_e$ agree at $x_0 = x_{2L}$,
the Wirtinger inequality{\footnote{The Wirtinger inequality is
just the usual Poincare inequality which bounds the $L^2$ norm in
terms of the $L^2$ norm of the derivative; i.e.,
$\int_{0}^{2\pi}f^2dt\leq 4\, \int_{0}^{2\pi}(f')^2dt$ provided
$f(0)=f(2\pi)=0$. }}  bounds $\dist^2 (\gamma_e ,
\tilde{\gamma}_e)$ in terms of
\begin{equation}    \label{e:tb}
    \int \, \left| (\tilde{\gamma}_e \circ P)' - \tilde{\gamma}_e' \right|^2 \leq
    2\, \int \,
\left| (\tilde{\gamma}_e' \circ P) \, P'  - \tilde{\gamma}_e'
\circ P \right|^2 +  2\, \int \, \left| \tilde{\gamma}_e' \circ P
- \tilde{\gamma}_e' \right|^2 \, .
\end{equation}
We will bound both terms on the right hand side of \eqr{e:tb} in
terms of $\int |P'-1|^2$ and then appeal to \eqr{e:ta}. To bound
the first term, use that $|\tilde{\gamma}_e'|$ is (a constant)
$\leq L$ to get
\begin{equation}
\int \, \left| (\tilde{\gamma}_e' \circ P) \, P'  -
\tilde{\gamma}_e' \circ P \right|^2 \leq   L^2 \int |P'-1|^2 \, .
\end{equation}
To bound  the second integral, we will use that when $x$ and $y$
are points in $\SS^1$ that are {\emph{not}} separated by a break
point, then $\tilde{\gamma}_e$ is a geodesic from $x$ to $y$ and,
thus, $\tilde{\gamma}_e''$ is normal to $M$ and
$|\tilde{\gamma}_e''| \leq |\tilde{\gamma}_e'|^2 \, \sup_M |A|
\leq \frac{L^2}{16}$. Therefore, integrating $\tilde{\gamma}_e''$
from $x$ to $y$ gives
\begin{equation}
\label{e:lip}
    |\tilde{\gamma}_e' (x) -
\tilde{\gamma}_e'(y)| \leq |x-y| \, \sup |\tilde{\gamma}_e''|
\leq \frac{L^2}{16} \, |x-y| \,
 .
\end{equation}
 Divide $\SS^1$ into two sets, $S_1$ and $S_2$, where $S_1$ is the set
of points within distance $(\pi \, \int |P'-1|^2 )^{1/2}$ of a
break point for $\tilde{\gamma}_e$. Since $P(x_0) = x_0$, arguing
as in \eqr{e:holder} gives $|P(x) - x| \leq (\pi \, \int |P'-1|^2
)^{1/2}$. Thus, if $x \in S_2$, then $\tilde{\gamma}_e$ is smooth
between $x$ and $P(x)$. Consequently, \eqr{e:lip} gives
\begin{equation}    \label{e:19}
    \int_{S_2} \,
\left| \tilde{\gamma}_e' \circ P  - \tilde{\gamma}_e' \right|^2
\leq  \frac{L^4}{256} \, \int_{S_2} \,   |P(s)-s|^2 \leq
\frac{L^4}{64} \, \int  \, |P'-1|^2 \, ,
\end{equation}
where the last inequality used the Wirtinger inequality. On the
other hand,
\begin{equation}    \label{e:110}
    \int_{S_1} \,
\left| \tilde{\gamma}_e' \circ P  - \tilde{\gamma}_e' \right|^2
\leq    4\, L^2 \, \Length (S_1)  \leq  8 \, L^3 \, \left(\pi \,
\int |P'-1|^2 \right)^{1/2} \, ,
\end{equation}
completing the proof of property (3).

\vskip1mm  We show (2) in Appendix \ref{s:aB}.

\vskip1mm  To prove property (4), we will argue by contradiction.
Suppose therefore that there exist $\epsilon > 0$ and a sequence
$\gamma_j \in \Lambda$ with  $\Energy (\Psi (\gamma_j)) \geq
\Energy (\gamma_j) - 1/j$ and $\dist (\gamma_j , G) \geq \epsilon
> 0$; note that the second condition implies a positive lower
bound for $\Energy (\gamma_j)$.   Observe next that the space
$\Lambda$ is compact{\footnote{Compactness of $\Lambda$ follows
 since  $\sigma \in \Lambda$ depends continuously on
the images of the $L$ break points in the compact manifold $M$.}}
and, thus, a
 subsequence of the $\gamma_j$'s must converge to some $\gamma \in
 \Lambda$.  Since property (3)
implies that $\dist (\gamma_j , \Psi (\gamma_j)) \to 0$, the $\Psi
(\gamma_j)$'s also converge to $\gamma$.  The continuity of
$\Psi$, i.e., property (2) of $\Psi$, then implies that $\Psi
(\gamma) = \gamma$.  However, this implies that $\gamma \in G$
since the only fixed points of $\Psi$ are immersed closed
geodesics. This last fact, which was used already by Birkhoff (see
section $2$ in \cite{Cr}), follows immediately from Lemma
\ref{l:sc} and \eqr{e:ta}.  However, this would contradict that
the $\gamma_j$'s remain a fixed distance from any such closed
immersed geodesic, completing the proof of (4).

\appendix

\section{Proof of Lemma \ref{l:sc}}  \label{s:B}

We will need a simple consequence of  (M1) and (M3)  in Section
\ref{s:existe}.

\begin{Lem}     \label{l:norm}
If $x,y \in M$, then $\left| (x-y)^{\perp} \right| \leq |x-y|^2$,
where $(x-y)^{\perp}$ is the normal component to $M$ at  $y$.
\end{Lem}

\begin{proof}
If $|x-y|\geq 1$, then the claim is clear. Assume therefore that
$|x-y| < 1$ and $\alpha:[0,\ell] \to M$ is a minimizing unit speed
geodesic from $y$ to $x$ with $\ell \leq 2 \,|x-y|$. Let $V$ be
the unit   normal vector $V=(x-y)^{\perp}/|(x-y)^{\perp}|$, so
$\langle \alpha'(0) , V \rangle = 0$, and observe that
\begin{align}
|(x-y)^{\perp}| &= \int_0^{\ell} \langle \alpha' (s) , V \rangle
\, ds = \int_0^{\ell} \,   \langle   \alpha'(0) + \int_0^s \,
\alpha'' (t) \, dt   \, , V \rangle  \, ds  \leq \int_0^{\ell} \,
\int_0^s \, \left| \alpha'' (t) \right| \, dt \,
ds  \notag \\
&\leq \int_0^{\ell} \, \int_0^s \, |A (\alpha (t))| \, dt  \, ds
\leq \frac{1}{2} \, \ell^2 \, \sup_M |A| \leq |x-y|^2 \, .
\end{align}
\end{proof}

\begin{proof}
(of Lemma \ref{l:sc}).
 Integrating by parts and using  that $\sigma_1$ and $\sigma_2$
are equal on $\partial I$ gives
\begin{equation}    \label{e:tma}
    \int_{I} |\sigma_1'|^2 - \int_{I} |\sigma_2'|^2 -
     \int_{I} \left| (\sigma_1 - \sigma_2)' \right|^2
    =   - 2 \, \int_{I}   \langle (\sigma_1 - \sigma_2)  ,  \sigma_2'' \rangle \equiv \kappa \, .
\end{equation}
The lemma will follow by bounding $|\kappa|$   by $\frac{1}{2} \,
        \int_{I} \left| (\sigma_1 - \sigma_2)' \right|^2$ and appealing to Wirtinger's inequality.

 Since $\sigma_2$ is a geodesic on $M$,
  $\sigma_2''$ is normal to $M$ and
 $| \sigma_2''| \leq  |\sigma_2'|^2 \, \sup_M |A| \leq \frac{|\sigma_2'|^2}{16}$.
 Thus,  Lemma \ref{l:norm} gives
\begin{equation}    \label{e:almosttan2a}
   \left| \langle (\sigma_1 - \sigma_2)  ,  \sigma_2'' \rangle \right|
   \leq |(\sigma_1 - \sigma_2)^{\perp}| \, \frac{|\sigma_2'|^2}{16}
   \leq
     |\sigma_1 - \sigma_2|^2 \, \frac{|\sigma_2'|^2}{16}
        \, .
\end{equation}
 Integrating \eqr{e:almosttan2a}, using that  $|\sigma_2'|$ is constant with $|\sigma_2'| \,
 \Length (I)  \leq 2 \pi$,  and applying
 Wirtinger's inequality gives
\begin{equation}    \label{e:gotpsia}
    \left| \kappa \right| \leq
    \frac{|\sigma_2'|^2}{8} \,
\int_{I} |\sigma_1 - \sigma_2|^2
 \leq \frac{|\sigma_2'|^2}{8} \,
     \left( \frac{ \Length (I) }{\pi} \right)^2 \, \int_{I}
                 |(\sigma_1 -\sigma_2)'|^2
                 \leq
 \frac{1}{2} \,
        \int_{I} \left| (\sigma_1 - \sigma_2)' \right|^2 \, .
\end{equation}
\end{proof}

\section{The continuity of $\Psi$}  \label{s:aB}

\begin{Lem} \label{l:cty}
Let $\gamma : \SS^1 \to M$  be a $W^{1,2}$ map  with $\Energy
(\gamma)   \leq L$.  If $\gamma_e$ and $\tilde{\gamma}_e$ are
given by applying steps $(A_1)$ and $(B_1)$ to $\gamma$, then the
map $\gamma \to \tilde{\gamma}_e$ is continuous from $W^{1,2}$ to
$\Lambda$ equipped with the $W^{1,2}$ norm.
\end{Lem}

\begin{proof}
It follows from \eqr{e:holder} and the energy bound that $\dist_M
(\gamma (x_{2j}) , \gamma ( x_{2j+2}) ) \leq 2\pi$ for each $j$
and thus we can apply step $(A_1)$.  The lemma will   follow
easily from two observations:
\begin{enumerate}
\item[(C1)]  Since $W^{1,2}$ close curves are also $C^0$ close
(cf. footnote $1$), it follows that the points $\gamma_e (x_{2j})
= \gamma (x_{2j})$ are continuous with respect to the $W^{1,2}$
norm.
 \item[(C2)] Define
$\Gamma \subset M \times M$ by $ \Gamma = \{ (x,y) \in M \times M
\, | \, \dist_M (x,y) \leq 4\pi \} \, $, and define a map
$H:\Gamma \to C^1([0,1],M)$ by letting $H(x,y): [0,1] \to M$ be
the linear map from $x$ to $y$. Then the map $H$ is continuous on
$\Gamma$. Furthermore, the map $t \to H(x,y)(t)$ has uniformly
bounded first and second derivatives $|\partial_t H(x,y)| \leq 4
\pi$ and $|\partial_t^2 H(x,y)| \leq   \pi^2$; the second
derivative bound comes from (M1).
\end{enumerate}
To prove the lemma, suppose   that
   $\gamma^1 $ and $\gamma^2 $ are non-constant curves in
   $\Lambda$
(continuity  at the constant maps is obvious).  For $i = 1,2$ and
$j=1 , \dots , L$, let $a^i_j$ be the distance in $M$ from
$\gamma^i (x_{2j})$ to $\gamma^i (x_{2j+2})$.  Let $S^i =
\frac{1}{2\pi} \, \sum_{j=1}^L a^i_j$ be the speed of
$\tilde{\gamma}^i_e$, so that $|(\tilde{\gamma}^i_e)'| = S^i$
except at the $L$ break points. By (C1),  the $a^i_j$'s  are
continuous functions of $\gamma^i$ and, thus, so are $S^1$ and
$S^2$.  Moreover, (C1) and (C2) imply that $\gamma_e^1$ and
$\gamma_e^2$ are $C^1$-close on each interval $[x_{2j},x_{2j+2}]$.
Thus, we have shown that $\gamma \to \gamma_e$ is continuous.

To show that $\gamma_e \to \tilde{\gamma}_e$ is also continuous,
we will show that the $\tilde{\gamma}^i_e$'s are close when the
$\gamma_e^i$'s are.   Since the point $x_0 = x_{2L}$ is fixed under the reparametrization,
this will follow from applying Wirtinger's inequality
to $(\tilde{\gamma}^1_e -
\tilde{\gamma}^2_e) -  ( \tilde{\gamma}^1_e -
\tilde{\gamma}^2_e)(x_0)$
once we show that $\int_{\SS^1} |(\tilde{\gamma}^1_e -
\tilde{\gamma}^2_e)'|^2$ can be made small.

The piecewise linear curve
$\tilde{\gamma}^i_e$ is linear on the intervals
\begin{equation} \label{e:where}
    I^i_j = \left[ \frac{1}{S^i} \, \sum_{\ell < j} a^i_{\ell} \, , \,
        \frac{1}{S^i} \, \sum_{\ell \leq j} a^i_{\ell} \right] \, .
   \end{equation}
   Set $I_j = I^1_j \cap I^2_j$.
   Observe first that
    since the intervals $I^i_j$ in \eqr{e:where} depend
   continuously on  $\gamma_e^i$,  the measure of the complement
   $\SS^1 \setminus \left[ \cup_{j=1}^L I_j \right]$ can be made small, so that
\begin{equation} \label{e:ij2}
    \int_{\SS^1 \setminus \left[ \cup I_j \right]} \, \, \left| (\tilde{\gamma}_e^1 -  \tilde{\gamma}_e^2)'
    \right|^2 \leq  4\,  L^2 \, \Length \, \left( \SS^1 \setminus \left[ \cup I_j
    \right] \right)
\end{equation}
can also be made small.  We will divide the $I_j$'s into two groups, depending on the size of $a^1_j$.
  Fix some $\epsilon > 0$ and suppose first that $a^1_j < \epsilon$; by continuity, we can assume that
  $a^2_j < 2\epsilon$.  For such  a $j$, we get
  \begin{equation}
    \int_{I_j}  \left| (\tilde{\gamma}_e^1 -  \tilde{\gamma}_e^2)'
    \right|^2 \leq 2 \, \int_{I_j^1}  \left| (\tilde{\gamma}_e^1)' \right|^2 + 2 \int_{I_j^2}
    \left| (  \tilde{\gamma}_e^2)' \right|^2 \leq 2 \, L \, \left( a^1_j + a^2_j \right) \leq 6 \, \epsilon \, L \, .
      \end{equation}
      Since there are at most $L$ breaks, summing over these intervals
      contributes at most $6\epsilon \, L^2$ to the energy of $(\tilde{\gamma}_e^1 -
      \tilde{\gamma}_e^2)$.

      The last case to consider is an $I_j$ with $a^1_j \geq \epsilon$; by continuity, we can assume that
  $a^2_j \geq \epsilon/2$.   In this case, $\tilde{\gamma}_e^i$
      can be written on $I_j$ as the composition
    $\gamma_e^i \circ
   P^i_j$ where $\left|
   (P^i_j)'\right| = 2 \pi \, S^i/ (L a^i_j)$.  Furthermore, $P^1_j$ and $P^2_j$ both map $I_j$ into $[x_{2j}, x_{2j+2}]$
   and
   \begin{equation} \label{e:ij}
    \int_{I_j} \left| (\tilde{\gamma}_e^1 -  \tilde{\gamma}_e^2)'
    \right|^2 =  \int_{I_j} \left| ({\gamma}_e^1 \circ P_j^1 -  {\gamma}_e^2 \circ
    P_j^2)' \right|^2 \, .
    \end{equation}
Finally, this can be made small since the speed $\left|
   (P^i_j)'\right|$ is continuous{\footnote{The speed is continuous because of the lower bound for the $a^i_j$'s.}}    in $\gamma^i$ and the ${\gamma}_e^i$'s are
   $C^2$ bounded and $C^1$ close on $[x_{2j}, x_{2j+2}]$.
   Therefore, the integral over these intervals can also be made
   small since there are at most $L$ of them.
\end{proof}

The next result shows that $\Psi$ preserves the homotopy class of
a sweepout.

\begin{Lem}     \label{l:homotopy}
Let $\gamma \in \Omega$ satisfy
\begin{equation}
    \max_t \, \, \Energy \, (\gamma (\cdot , t) ) \leq L \, .
\end{equation}
  If $\gamma_e$ and $\tilde{\gamma}_e$ are given by applying
steps $(A_1)$ and $(A_2)$ to each $\gamma (\cdot , t)$, then
$\gamma , \, \gamma_e$ and $\tilde{\gamma}_e$ are all homotopic in $\Omega$.
\end{Lem}

\begin{proof}
Given $x,y \in M$ with $\dist_M (x,y) \leq 4\pi$, let $H(x,y):
[0,1] \to M$ be the linear map from $x$ to $y$  as in (C2). It
follows that
\begin{equation}
    F(x,t,s) = H(\gamma (x,t) , \gamma_e (x,t)) (s)
\end{equation}
is an explicit homotopy with $F(\cdot , \cdot , 0) = \gamma$ and
 $F(\cdot , \cdot , 1) = \gamma_e$.

 For each $t$ with $\Length (\gamma_e (\cdot , t)) > 0$,
 $\gamma_e$ is given by $\gamma_e (\cdot , t) = \tilde{\gamma}_e (\cdot , t) \circ
 P_t$ where $P_t$ is a monotone reparametrization of $\SS^1$ that
 fixes $x_0 = x_{2L}$.
 Moreover, $P_t$ is continuous by \eqr{e:ta} and $P_t$ depends
 continuously on $t$ by Lemma \ref{l:cty}.  Since
 $x \to (1-s)P_t(x) + sx$ gives a homotopy from $P_t$ to the
 identity map on $\SS^1$, we conclude that
\begin{equation}
    G(x,t,s) =  \tilde{\gamma}_e  \, ( (1-s)P_t(x) + sx  , t)
\end{equation}
is an explicit homotopy with $G(\cdot , \cdot , 0) = \gamma_e$ and
 $G(\cdot , \cdot , 1) = \tilde{\gamma}_e$.  Note that $P_t$ is
 not defined when  $\Length (\gamma_e (\cdot , t))= 0$, but the
 homotopy $G$ is.
\end{proof}

\end{document}